\documentclass[letterpaper, 10pt]{article}
\usepackage{amsthm,amsmath,amsfonts}

\newtheorem{thm}{Theorem}[section]
\newtheorem{cor}[thm]{Corollary}
\newtheorem{lem}[thm]{Lemma}
\newtheorem{prop}[thm]{Proposition}
\newtheorem{clm}[thm]{Claim}
\theoremstyle{definition}

\theoremstyle{remark}

\newcommand{\Prob}[1]{\Pr\left[#1\right]}

\newcommand{\hG}{\widehat{G}}
\newcommand{\hV}{\widehat{V}}
\newcommand{\hE}{\widehat{E}}
\newcommand{\hX}{\widehat{X}}
\newcommand{\hU}{\widehat{U}}

\newcommand{\hW}{\widehat{W}}
\newcommand{\hw}{\widehat{w}}
\newcommand{\hn}{\widehat{n}}

\newcommand{\NhG}{N_{\widehat{G}}}
\newcommand{\Gnd}{\mathcal{G}_{n,d}}
\newcommand{\Gnp}{\mathcal{G}_{n,p}}
\newcommand{\avgd}{\widetilde{d}}

\title{Vertex Percolation on Expander Graphs}
\author{Sonny Ben-Shimon
\thanks{School of Computer Science, Raymond and Beverly Sackler
Faculty of Exact Sciences, Tel Aviv University, Tel Aviv 69978,
Israel. E-mail: sonny@post.tau.ac.il. Research conducted as part of
the author's Ph.D. thesis under the supervision of Prof. Michael
Krivelevich.}
\and Michael Krivelevich\thanks{School of Mathematical
Sciences, Raymond and Beverly Sackler Faculty of Exact Sciences, Tel
Aviv University, Tel Aviv 69978, Israel. E-mail:
krivelev@post.tau.ac.il. Research supported in part by USA-Israel
BSF Grant 2006-322, and by grant 526/05 from the Israel Science
Foundation.}}
\begin{document}
\maketitle
\begin{abstract}
We say that a graph $G=(V,E)$ on $n$ vertices is a $\beta$-expander
for some constant $\beta>0$ if every $U\subseteq V$ of cardinality
$|U|\leq \frac{n}{2}$ satisfies $|N_G(U)|\geq \beta|U|$ where
$N_G(U)$ denotes the neighborhood of $U$. In this work we explore
the process of deleting vertices of a $\beta$-expander independently
at random with probability $n^{-\alpha}$ for some constant
$\alpha>0$, and study the properties of the resulting graph. Our
main result states that as $n$ tends to infinity, the deletion
process performed on a $\beta$-expander graph of bounded degree will
result with high probability in a graph composed of a giant
component containing $n-o(n)$ vertices that is in itself an expander
graph, and constant size components. We proceed by applying the main
result to expander graphs with a positive spectral gap. In the
particular case of $(n,d,\lambda)$-graphs, that are such expanders,
we compute the values of $\alpha$, under additional constraints on
the graph, for which with high probability the resulting graph will
stay connected, or will be composed of a giant component and
isolated vertices. As a graph sampled from the uniform probability
space of $d$-regular graphs with high probability is an expander and
meets the additional constraints, this result strengthens a recent
result due to Greenhill, Holt and Wormald about vertex percolation
on random $d$-regular graphs. We conclude by showing that performing
the above described deletion process on graphs that expand
sub-linear sets by an unbounded expansion ratio, with high
probability results in a connected expander graph.
\end{abstract}

\section{Introduction}
In this paper we analyze the process of deleting vertices
independently at random from an expander graph and describe typical
properties and the structure of the resulting graph. We focus on the
case where the initial graph, $G$, is of bounded degree and the
deletion probability equals $n^{-\alpha}$, for any fixed $\alpha>0$,
where $n$ denotes the number of vertices in $G$. We are mainly
interested in investigating when the resulting graph with high
probability will possess some expansion properties as will be
discussed in Section \ref{ss:intromainresult}. In a recent paper of
Greenhill, Holt and Wormald \cite{GreHolWor2008}, the authors
perform a very similar analysis where the initial graph is sampled
from the uniform probability space of all $d$-regular graphs for
some fixed $d\geq 3$. Our current result, generalizing and improving
\cite{GreHolWor2008}, can be interpreted as providing sufficient
deterministic conditions on the initial graph that imply the result
of \cite{GreHolWor2008}. We are also able to prove some results when
the initial graph has an unbounded expansion ratio, and apply it to
the case of random $d$-regular graphs when $d=o(\sqrt n)$.
\subsection{Notation}
Given a graph $G=(V,E)$, the \emph{neighborhood} $N_G(U)$ of a
subset $U\subseteq V$ of vertices is the set of vertices defined by
$N_G(U)=\{v\notin U\;:\;v\hbox{ has a neighbor in }U\}$. For any
$f:\left[\lfloor\frac{n}{2}\rfloor\right]\rightarrow\mathbb{R}^+$,
we say that a graph $G=(V,E)$ on $n$ vertices is an
$f$-\emph{expander} if every $U\subseteq V$ of cardinality $|U|\leq
\frac{n}{2}$ satisfies: $|N_G(U)|\geq f(|U|)\cdot|U|$. When $f$ is a
constant function equal to some $\beta>0$ we say that $G$ is a
$\beta$-\emph{expander}. When a function $f:A\rightarrow
\mathbb{R}^+$ satisfies: $f(a)\geq c$ for any $a\in A$, where $c\geq
0$ is a constant, we simply write $f\geq c$.

Expanders in general are highly-connected sparse graphs. There are
many other notions and definitions of expander graphs different from
the one described above, some of which will be addressed in the
coming sections. Expander graphs is a subject of utmost importance
to the fields of both applied and theoretical Computer Science,
Combinatorics, Probability Theory etc. Monograph
\cite{HooLinWig2006} provides an excellent survey on expander graphs
and their applications.

In our setting, we start with a graph $G=(V,E)$ on $n$ vertices. We
delete every vertex of $V$ with probability $p=n^{-\alpha}$ for some
fixed  $\alpha>0$ independently at random. To simplify notation,
from here on, we will denote the resulting graph of this process by
$\hG=(\hV,\hE)$, and for every $X\subseteq V$, we denote by
$\hX=X\cap\hV$ the subset of $X$ that was not deleted by the
deletion process. We denote by $\hn$ the cardinality of $\hV$, by
$R$ the set of deleted vertices, i.e. $\hG=G[V\setminus R]$, and by
$r$ its cardinality, i.e. $\hn = n-r$. When considering the
neighborhood in the graph $\hG$ of a subset of vertices $U\subseteq
\hV$ we denote it by $\NhG(U)$.

The main research interest of this paper is the asymptotic behavior
of properties of the graph $\hG$ as we let the number of vertices,
$n$, grow to infinity. In this context, one needs to be precise when
formulating such claims. When stating an asymptotic claim for every
graph $G$ on $n$ vertices that satisfies a set of properties
$\mathcal{P}_n$ (the properties may depend on $n$), one actually
means that for every family of graphs $\mathcal{G}=\{G_n\}$, such
that $G_n$ is a graph on $n$ vertices satisfying $\mathcal{P}_n$,
there exists a value $n_0$ such that the claim is correct for every
$G_n$ where $n>n_0$. We say that an event $\mathcal{A}$ in our
probability space occurs with high probability (or w.h.p. for
brevity) if $\Prob{\mathcal{A}}\rightarrow 1$ as $n$ goes to
infinity. Therefore, from now on and throughout the rest of this
work, we will always assume $n$ to be large enough. We use the usual
asymptotic notation. That is, for two functions of $n$, $f(n)$ and
$g(n)$, we denote $f=O(g)$ if there exists a constant $C>0$ such
that $f(n)\leq C\cdot g(n)$ for large enough values of $n$; $f=o(g)$ or $f\ll g$
if $f/g\rightarrow 0$ as $n$ goes to infinity; $f=\Omega(g)$ if
$g=O(f)$; $f=\omega(g)$ or $f\gg g$ if $g=o(f)$; $f=\Theta(g)$ if
both $f=O(g)$ and $f=\Omega(g)$.

\subsection{Motivation}\label{ss:motivation}
Let $\Gnd$ denote the random graph model consisting of the uniform
distribution of all $d$-regular graphs on $n$ labeled vertices
(where $dn$ is even). One of the motivations of this paper is the
following result, recently proved by Greenhill, Holt and Wormald in
\cite{GreHolWor2008}.
\begin{thm}[Greenhill, Holt and Wormald \cite{GreHolWor2008}]\label{t:GreHolWor}
For every fixed $\alpha>0$ and fixed $d\geq 3$ there exists a
constant $\beta>0$, such that if $p=n^{-\alpha}$ and $G$ is a graph
sampled from $\Gnd$, then w.h.p. $\hG$ has a connected component of
size $n-o(n)$ that is a $\beta$-expander and all other components
are of bounded size. Moreover,
\begin{enumerate}
\item {if $\alpha>\frac{1}{2(d-1)}$, w.h.p. all small connected components of $\hG$ are isolated vertices.}\label{i:GHW2}
\item {if $\alpha\geq\frac{1}{d-1}$, w.h.p. $\hG$ is connected.}\label{i:GHW1}
\end{enumerate}
\end{thm}
Theorem \ref{t:GreHolWor} suggests a few questions that may be of
interest to address. First, one might consider the question whether
the deletion probability $p$ for which the desired properties hold
is best possible. This question has been answered in
\cite{GreHolWor2008}. Simple probabilistic arguments show that the
above result is indeed optimal in the sense that if we let
$\alpha=o(1)$, the largest component of $\hG$ will contain w.h.p.
many induced paths of length $O(1/\alpha)$, and hence cannot be an
expander. Next, one may ask what are the properties of random
$d$-regular graphs that make the above claim true. One of the
research motivations of this paper is precisely that, as will be
described below. Moreover, Item \ref{i:GHW1} of Theorem
\ref{t:GreHolWor} does not seem to be optimal due to the following
argument. As random $d$-regular graphs (for constant values of $d$)
w.h.p. locally look like trees (i.e. there are very few cycles of
constant length) it would seem natural to think that to disconnect
such a graph one would need to find the deletion probability that is
``just enough'' to disconnect a single vertex. A simple first moment
argument would imply that $\alpha>\frac{1}{d}$ should suffice. In
Section \ref{ss:ndlambda} we confirm this hypothesis in the more
general setting of pseudo-random $(n,d,\lambda)$-graphs. Lastly,
Theorem \ref{t:GreHolWor} does not consider the case of sampling a
random $d$-regular graph when $d=\omega(1)$, i.e. $d$ goes to
infinity with $n$. This setting is addressed in Section
\ref{s:unbounded}.

\subsection{Main result}\label{ss:intromainresult}
The main result of this paper states that the deletion of vertices
of an expander graph $G$ independently at random with probability
$n^{-\alpha}$ w.h.p. results in $\hG$ containing one giant
components that is in itself an expander graph. Moreover, the
expansion properties of $G$ imply a bound on the sizes of the small
connected components of $\hG$.

\begin{thm}\label{t:MainThm}
For every fixed $\alpha,c>0$ and fixed $\Delta>0$, there exists a
constant $\beta>0$, such that if $G$ is an $f$-expander graph on $n$
vertices of bounded maximum degree $\Delta$, and $f\geq c$, then
w.h.p. $\hG$ has a connected component of size $n-o(n)$ that is a
$\beta$-expander, and the rest of its connected components have at
most $K-1$ vertices, where
\begin{equation}\label{e:defnK}
K=\min\left\{u\;:\;\forall k\geq u\quad
kf(k)>\frac{1}{\alpha}\right\}.
\end{equation}
\end{thm}
We note that $K$ is well defined as $f\geq c$ implies that $K<
\frac{1}{c\alpha}$.

As mentioned in Section \ref{ss:motivation}, Theorem \ref{t:MainThm}
is optimal with respect to the deletion probability if we require
the giant component of the resulting graph to possess expansion
properties.

It is well known that for fixed $d\geq 3$ random $d$-regular graphs
are w.h.p. expander graphs. Thus our result strengthens Theorem
\ref{t:GreHolWor}, as will be formalized in Section
\ref{ss:RandRegGraphs}. It should be stressed that the techniques
used in the present paper and in \cite{GreHolWor2008} are quite
different. Whereas in \cite{GreHolWor2008} the analysis is done
directly in the so called Configuration Model in a probabilistic
setting, we rely upon a deterministic property of a graph, namely,
being an expander. The approach of first proving some claim under
deterministic assumptions, and then showing that these conditions
appear w.h.p. in some probability space, allows us to, arguably,
simplify the proof, and to get a strengthened result for families of
pseudo-random graphs and the random $d$-regular graph.

\subsection{Related work}\label{ss:relatedwork}
The process of random deletion of vertices of a graph received
rather limited attention, mainly in the context of faulty storage
(see e.g. \cite{Aloetal2002}), communication networks, and
distributed computing. For instance, the main motivation of
\cite{GreHolWor2008} is the SWAN peer-to-peer \cite{Holtetal2005}
network whose topology possess some properties of $d$-regular
graphs, and may have faulty nodes. Other works are mainly interested
in connectivity and routing in the resulting graph after performing
(possibly adversarial) vertex deletions on some prescribed graph
topologies.

The process of deleting edges, sometimes referred to by
\emph{edge-percolation} (or \emph{bond-percolation}) has been more extensively
studied. The main interest of edge-percolation is the existence of a
``giant component'', i.e. a connected component consisting of a
linear size of the vertices, in the resulting graph. When the
initial graph is taken to be $K_n$, edge-percolation becomes the
famous $\mathcal{G}(n,p)$ random graph model. In \cite{Aloetal2004,
FriKriMar2004, NachmiasPre} the edge percolation on an expander
graph is considered, the authors determine the threshold of the
deletion probability at which the giant component emerges w.h.p.. It
should be noted that in the context of this paper the expected
number of deleted vertices is far lower than permissable in order to
retain a giant component in the graph, as is clearly seen in Lemma
\ref{l:giantcomp}.

\subsection{Organization of paper}
The rest of the paper is organized as follows. In Section
\ref{s:ProofMainThm} we give a proof of Theorem \ref{t:MainThm}. We
proceed in Section \ref{s:Applications} to a straightforward
application of our result to expander graphs arising from
constraints on the spectrum of the graph. We continue in Section
\ref{ss:ndlambda} to the particular case of $(n,d,\lambda)$-graphs,
and under additional constraints on the graph compute the values of
$\alpha$ for which the resulting graph will w.h.p. stay connected or
will be composed of a giant component and isolated vertices. In
Section \ref{ss:RandRegGraphs} we show that a graph sampled from the
uniform probability space of $d$-regular graphs satisfies all
constraints, providing an alternative proof of the main result of
\cite{GreHolWor2008} and even improving it. As a final result, in
Section \ref{s:unbounded} we analyze the case of graphs of unbounded
expansion ratio for sub-linear sets with the same deletion
probability, and extend our result to random $d$-regular graphs
where $1\ll d\ll \sqrt n$. We conclude in Section \ref{s:conrem}
with a short summary and open problems for further research.

\section{Proof of Theorem \ref{t:MainThm}}\label{s:ProofMainThm}
Let $G$ be an $f$-expander, where $f\geq c$ for some constant $c>0$.
The number of deleted vertices, $r$, is clearly distributed by
$r\sim \mathcal{B}(n,p)$, and hence by the Chernoff bound (see e.g.
\cite{AlonSpencer2000}) $r$ is highly concentrated around its
expectation.
\begin{clm}\label{c:delcardinality}
W.h.p. $(1-o(1))n^{1-\alpha}\leq
r\leq(1+o(1))n^{1-\alpha}$.
\end{clm}

As for $\alpha> 1$ w.h.p. no vertices are deleted from the graph $G$ and
the proof of Theorem \ref{t:MainThm} becomes trivial, we will assume
from now on that $\alpha\leq 1$. Denote by $\hV_1,\ldots, \hV_s$ the
partition of $\hV$ to its connected components ordered in descending
order of cardinality. We call $\hV_1$ the \emph{big} component of
$\hG$, and $\hV_2,\ldots, \hV_s$ the \emph{small} components of
$\hG$.

\begin{lem}\label{l:giantcomp}
W.h.p. $|\hV_1|\geq(1-Cn^{-\alpha})n$ for any $C>\frac{1+c}{c}$.
\end{lem}

\begin{proof}
First, we show that $|\hV_1|>\frac{n}{2}$. Assume otherwise, and
take $\hW=\bigcup_{i=1}^j\hV_i$ for some $j\in[s]$ such that
$\frac{n}{4}\leq |\hW|\leq\frac{n}{2}$. Such a $j$ surely
exists. By our condition on $f$, we have that $|N_G(\hW)|\geq
c|\hW|=\Theta(n)$. But surely, $N_G(\hW)\subseteq R$, and hence, by
Claim \ref{c:delcardinality} $|N_G(\hW)|=o(n)$, a contradiction.
Now, set $\hU=\hV\setminus\hV_1$. From the above, it follows that
$|\hU|<\frac{n}{2}$. Clearly, $N_G(\hU)\subseteq R$, and
$|N_G(\hU)|\geq c|\hU|$. Putting these together yields that
$|\hU|\leq \frac{|R|}{c}$, and hence, by Claim
\ref{c:delcardinality}, w.h.p. $|V\setminus \hV_1|=|R\cup
\hU|\leq(1+o(1))\frac{1+c}{c}n^{1-\alpha}$, completing the proof.
\end{proof}

In a graph $H$, we call a subset of vertices $U\subseteq V(H)$
\emph{connected} if the corresponding spanning subgraph $H[U]$ is
connected. The following well known lemma (see e.g. \cite[Exercise
11, p.396]{Knuth69}) helps us to bound the number of connected
subsets of vertices in a graph of bounded maximum degree.

\begin{lem}\label{l:numconsubgraphs}
If $H=(V,E)$ is a graph of maximum degree $D$, then $V$ contains at
most $\frac{|V|(De)^k}{k}$ connected subsets of cardinality $k$.
\end{lem}

Keeping in mind that $\Delta$ is a constant, Lemma \ref{l:numconsubgraphs} turns out to be quite crucial to our
forthcoming calculations, for it allows us to bound probability of
events using union bound arguments by summing over connected
subgraphs of a prescribed cardinality instead of summing over all
subgraphs of the respective cardinality. We continue by showing that
w.h.p. all small connected components of $\hG$ must span less than
$K$ vertices.

\begin{lem}\label{l:smallconcomp}
W.h.p. every small connected component of $\hG$ is of cardinality at
most $K-1$.
\end{lem}
\begin{proof}
By Lemma \ref{l:giantcomp}, every small connected component of $\hG$
will be of cardinality at most $\bar{u}=O(n^{1-\alpha})$. Let
$U\subseteq V$ be a subset of vertices of cardinality $u\leq
\bar{u}$, and let $W=N_G(U)$ be its neighborhood in $G$, where
$|W|=w$. We first bound the probability that $\hU$ is a small
connected component of $\hG$ by the probability that all of $W$ was
deleted,
\begin{equation*}
\Prob{\exists j>1\hbox{ s.t. }\hU=\hV_j}\leq p^{w}\leq n^{-\alpha
f(u) u}.
\end{equation*}

By Lemma \ref{l:numconsubgraphs} we know that there are at most
$\frac{n(e\Delta)^u}{u}$ connected spanning subgraphs of cardinality
$u$, thus we can bound the probability of appearance of a small
connected component of cardinality $u$.
\begin{equation}\label{e:concompsizeu}
\Prob{\exists j>1\hbox{ s.t.
}|\hV_j|=u}\leq\frac{n(e\Delta)^u}{u}\cdot n^{-\alpha f(u) u}.
\end{equation}

Setting $\varepsilon=\min\{kf(k)-\frac{1}{\alpha}\;:\;k\geq K\}$,
the definition of $K$ implies $\varepsilon$ is a positive constant.
Applying \eqref{e:concompsizeu} and summing over all possible values
of $u$, we can bound the probability there will be in $\hG$ a small
connected component of cardinality at least $K$. First assume
$K<\lceil\frac{2}{c\alpha}\rceil$.
\begin{eqnarray*}\label{e:smallconcomp2}
\Prob{\exists j>1\hbox{ s.t. }|\hV_j|\geq K}
&\leq&\sum_{u=K}^{\bar{u}}\frac{n(e\Delta)^u}{u}\cdot n^{-\alpha f(u) u}\\
&\leq& \sum_{u=K}^{\lceil\frac{2}{c\alpha}\rceil-1}n^{1-\alpha uf(u)+o(1)} +\sum_{u=\lceil\frac{2}{c\alpha}\rceil}^{\bar{u}}n^{1-u\left(c\alpha-o(1)\right)}\\
&\leq&\left(\left\lceil\frac{2}{c\alpha}\right\rceil-K\right)n^{-\alpha\varepsilon+o(1)}+
O\left(n^{2-\alpha-\left\lceil\frac{2}{c\alpha}\right\rceil\left(c\alpha-o(1)\right)}\right)=o(1).
\end{eqnarray*}
Finally, if $K\geq \lceil\frac{2}{c\alpha}\rceil$ the above
computation becomes simpler as we are left with only the second
summand in the second line and the statement holds in that case as
well.
\end{proof}

Having shown that the small connected components of $\hG$ are w.h.p.
of bounded size, we move on to show that larger connected subsets of
$\hG$ w.h.p. expand.
\begin{lem}\label{l:consubgraphs}
W.h.p. every connected subset $U\subseteq \hV$ of $\hG$ of
cardinality $u$ s.t. $K\leq u\leq\frac{\hn}{2}$ satisfies
$|\NhG(U)|\geq \frac{\alpha c}{4} u$.
\end{lem}
\begin{proof}
Similarly to the proof of Lemma \ref{l:smallconcomp} let $W=N_G(U)$
and $\hW=\NhG(U)$ be the neighborhoods of $U$ in $G$ and $\hG$,
respectively, and let $w$ and $\hw$ denote their respective
cardinalities. By Lemma \ref{l:smallconcomp} we have that w.h.p.
every connected subset of vertices $U$ of cardinality $K\leq
u\leq\frac{4}{\alpha c}$ is not disconnected from the graph, and
thus has at least one edge leaving it. Setting $\eta=\frac{\alpha
c}{4}$ this implies that for every such connected subset $U$,
$\hw\geq \eta u$. Assuming $\frac{4}{\alpha c}< u\leq
\frac{\hn}{2}$, relying on $\hw\sim \mathcal{B}(w,1-p)$ we have
\begin{equation*}
\Prob{\hw<\eta u}\leq{w\choose{\lfloor\eta u\rfloor}}\cdot
p^{w-\lfloor\eta u\rfloor} \leq \left(\frac{ew}{\eta u}\right)^{\eta
u}p^{uf(u)-\eta u} \leq\left(\frac{e\Delta}{\eta}\right)^{\eta
u}n^{-\alpha u(f(u)-\eta)}.
\end{equation*}

To bound the probability there exists a connected subset in $\hG$ of
cardinality $u$ whose neighborhood contains less than $\eta u$
vertices, we apply the above with the union bound on all connected
subsets of $G$ from Lemma \ref{l:numconsubgraphs} as follows.
\begin{eqnarray*}
\frac{n(e\Delta)^u}{u}\cdot \left(\frac{e\Delta}{\eta}\right)^{\eta
u}\cdot n^{-\alpha u(f(u)-\eta)} &\leq& \eta^{-\eta u}\cdot
(e\Delta)^{u(\eta +1)}\cdot n^{1-\alpha u(f(u) -\eta)}\\
&\leq& \left(\frac{(e\Delta)^{\eta+1}}{n^{\alpha
c/4}\eta^{\eta}}\right)^u\cdot n^{-1}=o(n^{-1}).
\end{eqnarray*}
The inequality from the first to the second line follows from the
fact that $1-\alpha u(f(u)-\frac{\alpha c}{4})\leq  -(1+\frac{\alpha
u c}{4})$, or equivalently $\alpha u f(u) -\frac{\alpha u c}{4}
(1+\alpha)\geq \frac{\alpha u c}{2}\geq 2$ using that $\alpha\leq 1$
and $u> \frac{4}{\alpha c}$. Summing over all possible values of $u$
implies that w.h.p. there is no connected subset $U$ in $\hG$ of
cardinality at least $\frac{4}{\alpha c}$ that satisfies $|\NhG(U)|<
\eta|U|$, completing the proof.
\end{proof}
Lemma \ref{l:consubgraphs} states that w.h.p. all connected subsets of
$G[\hV_1]$ expand. As $G$ is of bounded maximum degree, this is
sufficient to imply that w.h.p. all subsets of $G[\hV_1]$ expand.

\begin{lem}\label{l:allsubgraphs}
W.h.p. $G[\hV_1]$ is a $\beta$-expander, where $\beta
=\frac{1}{\Delta}\cdot\min\{\frac{1}{K},\frac{\alpha c}{4}\}$.
\end{lem}
\begin{proof}
Set $\eta=\frac{\alpha c}{4}$ as defined in Lemma
\ref{l:consubgraphs}, and $\gamma=\min\{\frac{1}{K},\eta\}$. For
every $U\subseteq \hV_1$ of cardinality $|U|=u\leq K\leq 1/\gamma$,
trivially $|\NhG(U)|\geq 1\geq \gamma u$, as $U$ has at least one
edge emitting out of it. Assume $u>K$ and and denote by
$U_1,\ldots,U_t$ the decomposition of $U$ to its connected subsets,
and by $u_1,\ldots,u_t$ their respective cardinalities. As every
connected subset $U_i$ satisfies w.h.p. $|\NhG(U_i)|\geq \gamma u_i$
by Lemma \ref{l:consubgraphs}, it follows that w.h.p. $|\NhG(U)|\geq
\frac{\gamma}{\Delta}u$ completing the proof.
\end{proof}
Combining Lemmata \ref{l:giantcomp}, \ref{l:smallconcomp}, and \ref{l:allsubgraphs} completes the proof of Theorem \ref{t:MainThm}
\section{Applications to different expander graph families}\label{s:Applications}
\subsection{Expansion via the spectrum of a
graph}\label{ss:spectrum} The \emph{adjacency matrix} of a graph $G$
on $n$ vertices labeled by $\{1,\ldots, n\}$, is the $n\times n$
binary matrix, $A(G)$, where $A(G)_{ij}=1$ iff $i\sim j$. The
\emph{combinatorial Laplacian} of $G$ is the $n\times n$ matrix
$L(G)=D-A(G)$ where $D$ is the diagonal matrix defined by
$D_{i,i}=d_G(i)$. It is well known that for every graph $G$, the
matrix $L(G)$ is positive semi-definite (see e.g. \cite{Chung2004}),
and hence has an orthonormal basis of eigenvectors and all its
eigenvalues are non-negative. We denote the eigenvalues of $L(G)$ in
the ascending order by $0=\sigma_0\leq\sigma_1 \ldots\leq
\sigma_{n-1}$, where $\sigma_0$ corresponds to the eigenvector of
all ones. We denote by $\avgd=\avgd(G)$ the average degree of $G$,
and let $\theta=\theta(G)=\max\{|\avgd-\sigma_i|\;:\;i>0\}$. The
celebrated expander mixing lemma (see e.g. \cite{AlonSpencer2000})
and its generalization to the non-regular case (see e.g.
\cite{Chung2004}) state roughly that the smaller $\theta$ is, the
more random-like is the graph. This easily implies several
corollaries on the distribution of edges in the graph. In
particular, one can deduce the following expansion property of $G$
in terms of $\avgd$ and $\theta$. Its full proof can be found in
\cite{Chung2004}.
\begin{prop}\label{p:vertexboundarygen}
Let $G$ be an graph on $n$ vertices. Then $G$ is an
$h_{n,\avgd,\theta}$-expander, where
\begin{equation}\label{e:vertexboundarygen}
h_{n,\avgd,\theta}(i)=\frac{\avgd^2-\theta^2}{\theta^2+\avgd^2\frac{i}{n-i}}\;\hbox{
for }1\leq i\leq \left\lfloor\frac{n}{2}\right\rfloor.
\end{equation}
\end{prop}

Assume $G$ is a graph of bounded maximum degree, implying $\avgd
=O(1)$, and let $H(i)=i\cdot h_{n,\avgd,\theta}(i)$. Straightforward
analysis implies $H(i)$ is monotonically increasing for $1\leq i\leq
\left\lfloor \frac{n\theta}{\avgd+\theta}\right\rfloor$, and
monotonically decreasing for $\left\lceil
\frac{n\theta}{\avgd+\theta}\right\rceil\leq i\leq
\left\lfloor\frac{n}{2}\right\rfloor$. When
$\avgd-\theta>\varepsilon$ for some $\varepsilon>0$, we have that
$h_{n,d,\theta}\geq c$, where
$c=\frac{\avgd^2-\theta^2}{\avgd^2+\theta^2}$ is a constant
depending on $\avgd$ and $\theta$. We note that
$H(\left\lfloor\frac{n}{2}\right\rfloor)=O(n)\gg 1/\alpha$. Setting
$k=\frac{\theta^2}{(\avgd^2-\theta^2)\alpha}+1$, we have that $k\leq
\left\lfloor \frac{n\theta}{\avgd+\theta}\right\rfloor$ and
$$H(k)=\left(\frac{\theta^2}{(\avgd^2-\theta^2)\alpha}+1\right)\cdot\left(\frac{\avgd^2-\theta^2}{\theta^2+o(1)}\right)
>\frac{1}{\alpha}.$$ Our analysis of $H(i)$ implies that the value
$K$ defined in \eqref{e:defnK} satisfies $K\leq k$. Proposition
\ref{p:vertexboundarygen} thus enables us to apply Theorem
\ref{t:MainThm} to such graphs.
\begin{thm}\label{t:spectrumgen}
For every fixed $\alpha, \varepsilon>0$ and fixed $\Delta\geq 0$
there exists a constant $\beta>0$, such that if $G$ is a graph on
$n$ vertices of maximum degree $\Delta$, and
$\avgd-\theta>\varepsilon$, then w.h.p. $\hG$ has a connected
component of size $n-o(n)$ that is a $\beta$-expander, and all other
components are of cardinality at most
$\frac{\theta^2}{(\avgd^2-\theta^2)\alpha}$.
\end{thm}
\subsection{$(n,d,\lambda)$-graphs}\label{ss:ndlambda}
When the graph $G$ is $d$-regular, $L(G)=dI-A(G)$, and hence if
$\lambda_0\geq \lambda_1\geq\ldots\geq\lambda_{n-1}$ is the spectrum
of $A(G)$ we have that $\lambda_i=d-\sigma_i$. In the case of
$d$-regular graphs it is customary and, arguably, more natural to use
the spectrum of $A(G)$ rather than of $L(G)$ to address expansion
properties of the graph. As the largest eigenvalue of $A$ is clearly
$\lambda_0=d$ and it is maximal in absolute value, we have that
$\theta(G)=\max\{|\lambda_1(G)|,|\lambda_{n-1}(G)|\}$. In the case
of $d$-regular graphs it is customary to denote $\theta(G)$ by
$\lambda(G)=\lambda$, and to call such a graph $G$ an
$(n,d,\lambda)$-graph. For an extensive survey of fascinating
properties of $(n,d,\lambda)$-graphs the reader is referred to
\cite{KriSud2006}.

In the case of $(n,d,\lambda)$-graphs Proposition
\ref{p:vertexboundarygen} and Theorem \ref{t:spectrumgen} translate to the following.
\begin{prop}\label{p:vertexboundaryndlambda}
Let $G$ be an $(n,d,\lambda)$-graph, then $G$ is an
$h_{n,d,\lambda}$-expander, where
\begin{equation}\label{e:ndlambdaexp}
h_{n,d,\lambda}(1)=d;\qquad\hbox{and}\qquad
h_{n,d,\lambda}(i)=\frac{d^2-\lambda^2}{\lambda^2+d^2\frac{i}{n-i}}\;\hbox{
for }2\leq i\leq \left\lfloor\frac{n}{2}\right\rfloor.
\end{equation}
\end{prop}

\begin{thm}\label{t:ndlambda}
For every fixed $\alpha, \varepsilon>0$ and fixed $d\geq 3$ there
exists a constant $\beta>0$, such that if $G$ is an
$(n,d,\lambda)$-graph where $d-\lambda>\varepsilon$, then w.h.p.
$\hG$ has a connected component of size $n-o(n)$ that is a
$\beta$-expander, and all other components are of cardinality at
most $\frac{\lambda^2}{(d^2-\lambda^2)\alpha}$.
\end{thm}

The next two propositions allow us to get improved bounds on the
sizes of the small connected components of $\hG$. In Proposition
\ref{p:ndlambda2} we are interested in the values of $\alpha$ for
which $\hG$ is w.h.p. connected, and in Proposition
\ref{p:ndlambda3} in the values for which w.h.p. the small connected
components of $\hG$ are all isolated vertices. We compute these
values of $\alpha$ under some additional assumptions on the
$(n,d,\lambda)$-graph. Specifically, we require the graph to be
locally ``sparse'' and the \emph{spectral gap}, i.e. $d-\lambda$, to
be relatively large. Although these constraints may seem somewhat
artificial, they arise naturally in the setting of random
$d$-regular graphs as will be exposed in Section
\ref{ss:RandRegGraphs}. For any graph $G=(V,E)$ we denote by
\begin{equation*}
\rho(G,M)=\max\left\{\frac{e(U)}{|U|}\;:\;U\subseteq V \hbox{ s.t. }
|U|\leq M\right\},
\end{equation*}
where $e(U)$ denotes the number of edges of $G$ that have both
endpoints in $U$.

\begin{prop}\label{p:ndlambda2}
For every $\alpha>\frac{1}{d}$ and fixed $d\geq 3$, if $G$ is an
$(n,d,\lambda)$-graph satisfying $\lambda\leq
2\sqrt{d-1}+\frac{1}{40}$ and $\rho(G,d+29)\leq 1$, then w.h.p.
$\hG$ is connected.
\end{prop}
\begin{proof}
We prove that such a graph $G$ is an $f$-expander where $if(i)\geq
d$ for every $1\leq i\leq \lfloor\frac{n}{2}\rfloor$. Lemma
\ref{l:smallconcomp} will then imply that $K=1$, and hence $\hG$ is
w.h.p. connected. Proposition \ref{p:vertexboundaryndlambda}
guarantees that $f(i)\geq h_{n,d,\lambda}(i)$ (with
$h_{n,d,\lambda}$ as defined in \eqref{e:ndlambdaexp}). Taking
$i\geq 30$ and plugging our assumption on $\lambda$ in the
definition of $h_{n,d,\lambda}$, we have that for every $d\geq 3$
\begin{equation*}
f(30)\geq\frac{d^2-4(d-1)-\frac{\sqrt{d-1}}{10}-\frac{1}{1600}}{4(d-1)+
\frac{\sqrt{d-1}}{10}+\frac{1}{1600}+d^2\frac{30}{n-30}}>\frac{d}{30}.
\end{equation*}
The analysis of $H$ therefore guarantees that $if(i)\geq d$ for all
$i\geq 30$. Now, let $U$ be a subset of vertices of cardinality
$u\leq 29$, and set $s=|N_G(U)|$ and $w=|U\cup N_G(U)|$. It now
suffices to show that $s\geq d$ for such a set $U$. If $u=1$,
trivially $s=d$, as every vertex has $d$ neighbors. Now, our
assumption on $G$ implies that all triangles in the graph must be
edge disjoint. Taking $u=2$, if the two vertices in $U$ are
non-adjacent trivially $s\geq d$, and if they are adjacent, they
must have at most one common neighbor, implying $s\geq 2d-3\geq d$.
Taking $3\leq u\leq 29$, if $w\geq u+d$ we are done. Otherwise, the assumption on $G$ implies $e(U)\leq u$
and $e(U\cup N_G(U))\leq w$, and hence $du-u\leq du-e(U)\leq e(U\cup
N_G(U))\leq w=u+s$. This implies $s\geq u(d-2)\geq d$ which
completes the proof.
\end{proof}

For any graph $G$ we denote by $t(G)$ the number of triangles in
$G$. To analyze the values of $\alpha$ for which w.h.p. all small
connected components of $\hG$ are isolated vertices, we additionally
require that the number of triangles in $G$ is bounded by a certain
positive power of $n$. This requirement as well is quite natural in
the case of random $d$-regular graphs.

\begin{prop}\label{p:ndlambda3}
For every $\alpha>\frac{1}{2(d-1)}$ and fixed $d\geq 3$, if $G$ is
an $(n,d,\lambda)$-graph satisfying $\lambda\leq
2\sqrt{d-1}+\frac{1}{40}$, $\rho(G,39+2(d-1))\leq 1$, and
$t(G)=O(n^\frac{2d-3}{2(d-1)})$, then w.h.p. all small connected
components of $\hG$ are isolated vertices.
\end{prop}
\begin{proof}
Following the spirit of the proof of Proposition \ref{p:ndlambda2},
we would like to show that $G$ is an $f$-expander where $if(i)\geq
2(d-1)$ for every $2\leq i\leq \lfloor\frac{n}{2}\rfloor$, which
completes the proof by using Lemma \ref{l:smallconcomp}. By
Proposition \ref{p:vertexboundaryndlambda} we can assume that
$f(i)\geq h_{n,d,\lambda}(i)$. Our assumption on $\lambda$ gives
\begin{equation*}
f(40)\geq\frac{d^2-4(d-1)-\frac{\sqrt{d-1}}{10}-\frac{1}{1600}}{4(d-1)+
\frac{\sqrt{d-1}}{10}+\frac{1}{1600}+d^2\frac{40}{n-40}}>\frac{d-1}{20},
\end{equation*}
which in turn, using our analysis of $H$, implies $if(i)\geq 2(d-1)$
for $i\geq 40$. When trying to complete the proof by showing that
$f(i)\geq\frac{2d-1}{i}$ for $2\leq i\leq 39$, it turns out that in
this setting this will not be the case, as there can be small
subsets that violate this strict expansion requirement. Fortunately,
we can prove that there cannot be too many such subsets, which
allows us to prove the above probabilistic statement.

Let $U$ be a subset of vertices of cardinality $2\leq u\leq 39$, and
set $s=|N_G(U)|$ and $w=|U\cup N_G(U)|$. We call $U$
\emph{exceptional} if $s< 2(d-1)$. Let $x_i$ denote the number of
exceptional sets of cardinality $i$.

If $4\leq u\leq 39$, our assumption on $G$ implies $e(U)\leq u$ and
$e(U\cup N_G(U))\leq w$. It follows that $du-u\leq e(U\cup
N_G(U))\leq w=u+s$, implying $s\geq u(d-2)\geq 2(d-1)$. If $u=2$ and
the two vertices of $U$ are non-adjacent then $2d\leq e(U\cup
N_G(U))\leq w=2+s$, hence $s\geq 2(d-1)$. If the two vertices are
adjacent, but are not part of a triangle, then again $s\geq 2(d-1)$.
For $u=3$, if $U$ spans at most one edge, then $3d-1\leq e(U\cup
N_G(U))\leq w=3+s$, yielding $s\geq 2(d-1)$. If $U$ spans a triangle
and $d\geq 4$, as all triangles of $G$ must be edge disjoint we get
that $s\geq 3(d-2)\geq 2(d-1)$. If $U$ spans exactly two edges, easy
case analysis, relying on the fact that no small subgraph spans more
edges then vertices, shows that $s\geq 3d-5\geq 2(d-1)$.

Lemma \ref{l:smallconcomp} and the previous computation assure that
for $d\geq 4$ w.h.p. all small connected components have at most two
vertices, and for $d=3$, w.h.p. all small connected components have
at most three vertices.

We conclude by showing that since in both cases there are only a small number of exceptional
sets, w.h.p. all small connected components will be isolated
vertices. Similarly to the proof of Lemma \ref{l:smallconcomp}, we bound the
probability of appearance of a connected component of cardinality
$2$. The exceptional sets of cardinality $2$ are edges that
participate in a triangle and all triangles in $G$ are edge
disjoint, therefore there are $x_2=3t(G)$ such exceptional sets, and each has exactly $2d-3$ neighbors.
Going over all connected sets of $G$ of cardinality $2$, i.e. the edges of
$G$, we bound the probability that one of these sets becomes
disconnected.
\begin{eqnarray*}
\Prob{\exists j>1\hbox{ s.t.
}|\hV_j|=2}&\leq&x_2p^{2d-3}+\left(\frac{dn}{2}-x_2\right)p^{2(d-1)}\\
&\leq&
O\left(n^{\frac{2d-3}{2(d-1)}-\alpha(2d-3)}\right)+O\left(n^{1-2\alpha(d-1)}\right)=o(1).
\end{eqnarray*}
The above completes the proof when $d\geq4$. We are left with the
case of exceptional triples that may exist when $d=3$. Since the
exceptional sets of cardinality $3$ are the triangles in $G$, there
are exactly $x_3=t(G)$ such exceptional sets each having exactly $3$
neighbors. Very similarly to the preceding computation, we go over
all connected sets of $G$ of cardinality $3$, i.e. sets that span
two or three edges, and compute the probability that one of these
sets become disconnected. Recall that for $d=3$ we have that
$\alpha>\frac{1}{4}$.
\begin{eqnarray*}
\Prob{\exists j>1\hbox{ s.t.
}|\hV_j|=3}&\leq&x_3p^{3}+\left(3n-3t(G)\right)p^{4}\\
&\leq&
O\left(n^{\frac{3}{4}-3\alpha}\right)+O\left(n^{1-4\alpha}\right)=o(1).
\end{eqnarray*}
\end{proof}
Propositions \ref{p:ndlambda2} and \ref{p:ndlambda3} are easily seen
to be optimal in some sense, for if $\alpha\leq\frac{1}{d}$ or
$\alpha\leq \frac{1}{2(d-1)}$, then the expected number of isolated
vertices or edges respectively is greater than 1.

\subsection{Random $d$-regular graphs}\label{ss:RandRegGraphs}
Consider the random graph model consisting of the uniform
distribution on all $d$-regular graphs on $n$ vertices (where $dn$
is even), and denote this probability space by $\Gnd$. Assume
throughout this section that $d\geq 3$ is a constant. Let $G$ be a
graph sampled from $\Gnd$. Note that the multiplicity of the
eigenvalue $d$ of the graph $G$ is w.h.p. $1$ as $G$ is w.h.p.
connected and non-bipartite (see e.g. \cite{Wormald99}), hence
w.h.p. $\lambda(G)<d$. Friedman, confirming a conjecture of Alon,
gives an accurate evaluation of $\lambda(G)$ for most random
$d$-regular graphs when $d$ is a constant.
\begin{thm}[Friedman \cite{Fri2004}]\label{t:Fri2004}
For any $\varepsilon>0$ and fixed $d\geq 3$, if $G$ is sampled from $\Gnd$ then w.h.p.
\begin{equation}\label{e:rhobound}
\lambda(G)\leq 2\sqrt{d-1}+\varepsilon.
\end{equation}
\end{thm}
Combining Theorems \ref{t:ndlambda} and \ref{t:Fri2004}, implies
explicitly the first part of Theorem \ref{t:GreHolWor}.
\begin{cor}\label{c:GreHolWor1}
For every fixed $\alpha>0$ and fixed $d\geq 3$ there exists a
constant $\beta>0$, such that if $G$ is a graph sampled from $\Gnd$,
then w.h.p. $\hG$ has a connected component of size $n-o(n)$ that is
a $\beta$-expander and all other components are of cardinality at
most $\frac{4(d-1)}{\alpha(d-2)^2}+1$.
\end{cor}
The second part of Theorem \ref{t:GreHolWor} analyzes the values of
$\alpha$ for which w.h.p. the graph $\hG$ is  connected, and the
values of $\alpha$ for which the w.h.p. small connected components
are all isolated vertices. Plugging Theorem \ref{t:Fri2004} into
Theorem \ref{t:ndlambda}, as above, implies a similar result, but
not as strong.

To get Theorem \ref{t:GreHolWor} in full, and even to improve it, we
use Propositions \ref{p:ndlambda2} and \ref{p:ndlambda3}. To do so,
we state the following well known asymptotic properties of $\Gnd$
(see e.g. \cite{Wormald99}). Let $G$ be a graph sampled from $\Gnd$,
for any fixed $d\geq3$, then w.h.p. the minimal distance between two
cycles of constant length in $G$ is $\omega(1)$. This statement is
equivalent to saying that for every constant $M>1$ w.h.p.
$\rho(\Gnd,M)\leq 1$. Moreover, as $n$ tends to infinity
$t(\Gnd)\sim Poisson\left(\frac{(d-1)^3}{6}\right)$, and by so
Markov's inequality w.h.p. $t(\Gnd)=O(n^{\frac {2d-3}{2(d-1)}})$
(with room to spare). Now, using Propositions \ref{p:ndlambda2} and
\ref{p:ndlambda3} combined with Corollary \ref{c:GreHolWor1} we get
the desired result for $\Gnd$.
\begin{thm}\label{t:improvedGreHolWor}
For every fixed $\alpha>0$ and $d\geq 3$ there exists a constant
$\beta>0$, such that if $p=n^{-\alpha}$ and $G$ is a graph sampled
from $\Gnd$, then w.h.p. $\hG$ has a connected component of size
$n-o(n)$ that is a $\beta$-expander and all other components are of
cardinality at most $\frac{4(d-1)}{\alpha(d-2)^2}+1$. Moreover,
\begin{enumerate}
\item {if $\alpha>\frac{1}{2(d-1)}$, w.h.p. all small connected components of $\hG$ are isolated vertices.}\label{i:improvedGHW2}
\item {if $\alpha>\frac{1}{d}$, w.h.p. $\hG$ is connected.}\label{i:improvedGHW1}
\end{enumerate}
\end{thm}
It should be noted that Theorem \ref{t:improvedGreHolWor} improves
upon Theorem \ref{t:GreHolWor} for the values of $\alpha$
guaranteeing that $\hG$ stays connected w.h.p.. As mentioned in the
Section \ref{ss:ndlambda}, this improvement is best possible, for if
$\alpha\leq\frac{1}{d}$, then the expected number of isolated
vertices will be at least one, and by some standard concentration
arguments it can also be shown that the number of isolated vertices
is highly concentrated around this expectation. Hence, for
$\alpha\leq \frac{1}{d}$ the graph $\hG$ has isolated vertices, and is thus
disconnected, with some probability bounded away from $0$. As a
final note, it should be mentioned that in the original statement of
the main result of \cite{GreHolWor2008}, it is proved that w.h.p.
all small connected components are trees, and that for $\alpha>
\frac{1}{2(d-1)}$ w.h.p. the number of isolated vertices is
$o(n^{(d-2)/2(d-1)})$. These results as well can be derived from
simple probabilistic arguments based on properties of $\Gnd$, but we
omit these technical details.
\section{Unbounded expansion of small sets}\label{s:unbounded}
So far we have considered graphs of bounded maximum degree (and in
particular $d$-regular graphs for $d=O(1)$) that expand by a
constant factor. When considering graphs that expand sets of
sub-linear cardinality by an $\omega(1)$ factor (in particular in
such graphs $\delta(G)=\omega(1)$, i.e. the minimal degree of $G$
goes to infinity with $n$) a simple union bound argument implies the
following result. The proof is quite similar to those we have
previously presented, only in this case we can use a union bound
over all subsets of vertices with no need to go over all connected
subsets first, i.e. we do not make use of Lemma
\ref{l:numconsubgraphs} .
\begin{thm}\label{t:unboundedgraphs}
For every fixed $\alpha,c,\varepsilon>0$ if $G$ is an $f$-expander
graph on $n$ vertices where $f(u)=\omega (1)$ for every $u=o(n)$,
and $f\geq c$, then w.h.p. $\hG$ is a $(c-\varepsilon)$-expander.
\end{thm}
\begin{proof}
Let $U\subseteq V$ be a subset of vertices of cardinality
$u\leq\frac{n}{2}$, and let $W=N_G(U)$ be its neighborhood in $G$,
where $|W|=w$. Set $\beta=c-\varepsilon$, and let us denote a subset
of vertices $U$ as \emph{bad} if $\widehat{U}=U$ and
$\widehat{w}<\beta u$. If $\hG$ is not a $\beta$-expander then it
must contain such a bad set. We bound the probability of a subset
$U$ to be bad by
\begin{equation*}
\Prob{U\hbox{ is bad}}\leq{w\choose{\lfloor\beta u\rfloor}}\cdot
p^{w-\lfloor\beta u\rfloor}\leq\left(\frac{ew}{\beta
u}\right)^{\beta u}p^{u(f(u)-\beta)}.
\end{equation*}
Assuming $u=o(n)$, we have
\begin{equation*}
\Prob{\exists U\subseteq V\hbox{ s.t. }|U|=u\hbox{ and }U\hbox{ is
bad}}\leq{n\choose u}\left(\frac{ew}{\beta u}\right)^{\beta
u}p^{u(f(u)-\beta)} \leq
n^{u(1+\beta(1+\alpha)+o(1)-\alpha\omega(1))}=o(n^{-1}).
\end{equation*}
In the case that $\Theta(n)=u\leq\frac{n}{2}$, we have
\begin{equation*}
\Prob{\exists U\subseteq V\hbox{ s.t. }|U|=u\hbox{ and }U\hbox{ is
bad}} \leq {n\choose u}\left(\frac{ew}{\beta u}\right)^{\beta
u}p^{u(f(u)-\beta)} \leq n^{u(o(1)+o(1)-\alpha(c-\beta))}=o(n^{-1}).
\end{equation*}
Applying the union bound over all possible values of $u$ completes
the proof.
\end{proof}

It should be noted that Theorem \ref{t:unboundedgraphs} implies that
when $p=n^{-\alpha}$ for any fixed $\alpha>0$, $\hG$ is w.h.p. an
expander, and in particular stays connected as opposed to the case
of bounded maximum degree.

When $d=o(\sqrt n)$, Broder et al. \cite[Lemma 18]{BroEtAl99}
provide an upper bound on the second eigenvalue of most of the
$d$-regular graphs.
\begin{thm}[Broder et al. \cite{BroEtAl99}]\label{t:BroEtAl99}
For $d=o(\sqrt n)$, if $G$ is sampled from $\Gnd$ then w.h.p.
\begin{equation}\label{e:rhobound2}
\lambda(G)= O(\sqrt d).
\end{equation}
\end{thm}

Plugging Theorem \ref{t:BroEtAl99} into Proposition
\ref{p:vertexboundaryndlambda} assures that w.h.p. all conditions
needed in Theorem \ref{t:unboundedgraphs} are met when the graph
sampled from $\Gnd$ for $1\ll d\ll \sqrt n$, and hence we get the
following result.
\begin{thm}\label{t:Gndunbounded}
For every fixed $\alpha>0$ and $1\ll d\ll \sqrt n$ there exists a
constant $\beta>0$, such that if $G$ is a graph sampled from $\Gnd$,
then w.h.p. $\hG$ is a $\beta$-expander.
\end{thm}
When sampling a graph from the binomial random graph model $\Gnp$
(i.e. the probability space of all graphs on $n$ labeled vertices,
where each pair of vertices is chosen to be an edge independently
with probability $p$) with $p=\frac{d}{n}$ for $d=\Omega(\sqrt n)$,
the graph is easily seen to be ``almost $d$-regular'' as all degrees
of the vertices are highly concentrated around $d$. Furthermore, it
can be easily shown that when the initial graph is sampled from
$\Gnp$ with the prescribed values of $p$, a similar claim to Theorem
\ref{t:unboundedgraphs} holds. Therefore, one should expect Theorem
\ref{t:unboundedgraphs} to extend to values of $d=\Omega(\sqrt n)$,
but unfortunately, the techniques that are commonly used to deal
with random regular graphs seem to fail for these higher values of
$d$.

We note that in \cite{BSKriPre} the authors prove a result on the
distribution of edges in $\Gnd$ for $d=o(\sqrt n)$, that can be
easily used to derive vertex-expansion properties of $\Gnd$ for
$1\ll d\ll \sqrt n$, and combined with Theorem
\ref{t:unboundedgraphs} provides an alternative proof of Theorem
\ref{t:Gndunbounded}.

\section{Concluding remarks and open problems}\label{s:conrem}
In this paper we analyzed the process of deleting uniformly at
random vertices from an expander graph. We have shown that for small
enough deletion probabilities the resulting graph w.h.p. retains
some expansion properties (if not in the graph itself then in its
largest connected component). We have also proved that for these
deletion probabilities w.h.p. all small connected components must be
of bounded size. Lastly, we have shown how this result can be
applied to the random $d$-regular graph model for $d=o(\sqrt n)$.

In Section \ref{ss:RandRegGraphs}, in order to apply our results
from previous sections to the case of random $d$-regular graphs, we
made use of several theorems that describe some properties that
occur w.h.p. in graphs that are sampled from $\Gnd$, such as Theorem
\ref{t:Fri2004} of Friedman \cite{Fri2004}. This very strong result,
whose proof is far from simple, seems to be an overkill to prove our
claims. One could go about by showing that graphs from $\Gnd$ w.h.p.
possess some expansion property (by analyzing the model directly
using, e.g., the Configuration Model or the Switching Technique) and
then by applying Theorem \ref{t:MainThm} directly. This method would
undoubtedly provide a proof that does not require any ``heavy duty
machinery'', but does require more meticulous computations.
Nonetheless, we hope that the reader finds the use of the connection
between spectral graph theory and expansion properties (or
pseudo-randomness of a graph) to be both elegant and concise.

In light of Theorem \ref{t:Gndunbounded} it would be interesting to
analyze the expansion properties of random $d$-regular graphs for
$d=\omega(1)$ for higher values of $p$, i.e. taking $p=n^{-o(1)}$ , as for $d=\omega(1)$ it is
no longer true trivially that for these values of $p$ w.h.p. there will
be long induced paths in $\hG$.

\subsection*{Acknowledgements}
The authors would like to thank the anonymous referees for their
helpful corrections and comments. The first author would like to
thank Itai Benjamini for discussing this problem with him.

\end{document}